\documentclass[leqno]{article}
\usepackage{latexsym,amsmath,amsfonts,mathrsfs,amssymb,amsthm}
\usepackage{yhmath}
\usepackage[usenames]{color}
\usepackage{graphicx}
\def\R{\mathbb R}
\def\N{\mathbb N}

\def\Z{\mathbb Z}

\def\a{\alpha}
\def\e{\epsilon}
\def\d{\delta}

\def\H{{\cal H} }
\def\K{{\cal K} }
\def\A{{\cal A} }
\def\F{{\cal F} }
\def\be{\begin{equation}}
\def\ee{\end{equation}}
\def\bs{\backslash}
\def\qed{\hfill$\Box$\bigskip}
\def\nd{\noindent Proof. }
\numberwithin{equation}{section}
\newtheorem{lem}[equation]{Lemma}
\newtheorem{pro}[equation]{Proposition}
\newtheorem{defn}[equation]{Definition}
\newtheorem{thm}[equation]{Theorem}
\newtheorem{cor}[equation]{Corollary}
\newtheorem{rem}[equation]{Remark}
\begin{document}
\bigskip

\begin{center}{\Large \textbf{Limits of topological minimal sets with finitely generated coefficient groups}}
\end{center}
\bigskip

\centerline{\large Xiangyu Liang}

\bigskip

\begin{center}Institut Camille Jordan, Universit\'e Claude Bernard Lyon 1, 43 boulevard du 11 novembre 1918, 69622 Villeurbanne Cedex, France\end{center}
\centerline{xiangyuliang@gmail.com}

\vskip 1cm

\centerline {\large\textbf{Abstract.}}We prove that the (local) Hausdorff limit of topological minimal sets (with finitely generated coefficient group) are topologically minimal. The key idea is to reduce the homology group on the space to the homology group on the sphere, and then reduce the homology group on the sphere to a finitely representable one, by "glueing" grids with small measure to block local elements in the homology group.

\bigskip


\bigskip


\section{Introduction}

It is frequently asked that given a sequence of sets, measures or functions that admit a certain important property, whether this property is kept while passing to the limit. For instance, this kind of compactness help to prove existence results for many minimizing problems in geometry measure theory. Also, in various models in geometric measure theory, a typical method to study the local structure around a point $x$ of a set $E$ is to look at the "blow-up" limits at the point, which are limits of subsequences of 
\be E_r=\frac 1r (E-x), r\to 0.\ee
This is also similar to the tangent measure in a measure theoretical setting. 

A blow-up limit of $E$ at $x$ describes the asymptotic behavior of $E$ around $x$ at small scales. The study of blow-up limits for sets, as well as functions, measures, etc., is widely used in the study of regularity and classification of singularities for many problems in geometric measure theory, geometric analysis, and calculus of variations. See works of Besicovitch, David, Hamilton, Kenig, Mattila, Preiss, Simon, Toro,  etc.. In most cases, dilatations always keep useful properties of the set. However, to carry on the study for limit behavior, one must prove that these properties pass to the limit. 

In this article, we discuss this problem for topological minimality of sets. The notion of minimal sets was initially introduced by Almgren \cite{Al76} to study soap films (or Plateau's problem in general dimensions and codimentions) in a setting of sets. Plateau's problems aims at understanding existence and regularity for physical objects that admitting certain minimizing properties, which is one of the central interests in geometric measure theory.

In Almgren's definition, a closed set $E$ is $d-$dimensional Almgren-minimal when there is no deformation $F=\varphi(E)$, where $\varphi$ is Lipschitz and $\varphi(x)-x$ is compactly supported, for which the Hausdorff measure $\H^d(F)$ is smaller than $\H^d(E)$. See Definition 2.9 for the precise definition.

The idea of minimizing measure among deformations corresponds with physical intuition for the formation of soap films. On the other hand, deformation is not an extrinsic property for sets---one set can have many different parametrizations. This brings many mathematical obstacles for proving results that seem to be obvious in physics. For example, we do not have any good existence result for Almgren minimal sets.

Another slightly stronger notion of minimal sets is the notion of topological minimal sets (introduced by the author in \cite{topo}). It is also in the setting of sets, but instead of minimizing Hausdorff measure among compact deformations, one asks that a topological minimal set admits a minimal measure among all sets that keep some topological property of the set. A simplified version for topological minimal sets of dimension $d$ (in $\R^n$, with coefficient group $G=\Z$) is the following:

\begin{defn}Let $U\subset\R^n$ be open. Let $E\subset U$ be relatively closed and has locally finite $d-$dimensional Hausdorff measure. Then $E$ is said to be $d-$dimensional topologically minimal in $U$ if 

\be \H^d(E\bs F)\le \H^d(F\bs E)\ee

for each closed set $F\subset U$ such that there exists a compact ball $B\subset U$ with the following properties:

$1^\circ$ $E\bs B=F\bs B$;

$2^\circ$ For each $n-d-1$-simplicial cycle $\gamma\subset U\bs(E\cup B)$, if $\gamma$ represents a non-zero element in $H_{n-d-1}(U\bs E,\Z)$, then it also represents a non zero element in $H_{n-d-1}(U\bs F,\Z)$.

Such a $F$ is called a topological competitor of dimension $d$ for $E$ in $B$.
\end{defn}
 
 A more general definition will be given in Definitions 2.14 and 2.18. When $d=n-1$, this is the Mumford-Shah minimal set defined in \cite{DJT}.
 
This definition might be physically less intuitive than that of Almgren minimal sets. However, for topological minimal sets, one can prove many good properties which we do not know how to prove for Almgren's minimal sets. For instance, existence results (cf. \cite{topo} Theorems 4.2 and 4.28), and the topological minimality of the product of a topological minimal set with $\R^n$ (cf. \cite{topo} Proposition 3.23). In addition, for many known Almgren minimal sets, their Almgren minimality was in fact proved by proving this stronger topological property (cf. e.g. \cite{LM94},\cite{Br91},\cite{YXY}). 
 
 Hence it would be also interesting to study local behaviour for topological minimal sets, in particular, the blow-up limits of such sets. As a first step, we will prove the following theorem:
 
\begin{thm}Let $E$ be the Hausdorff limit of a sequence of $d-$dimensional topological minimal sets $E_k,k\in \N$ in an open set $U\subset\R^n$ . Then $E$ is topologically minimal of dimension $d$ in $U$.\end{thm}

See Theorem 3.1 (the main theorem) for a more general version.

As a direct corollary, we will prove that the blow-up limits of a topological minimal set are all topological minimal cones. Also, the theorem makes it possible to use compactness argument in many circumstances. 

The idea of the proof of the main theorem is the following.

Suppose $E$ is not topologically minimal. Thus there exists a ball $B$ and a competitor $F$ of $E$ in $B$, such that $\H^d(E\cap B)>\H^d(F\cap B)$. We want to use $F$ to construct better competitors $E_k$ for $k$ large. A natural idea is to glue $E_k\bs B$ and $F\cap B$ together. 

Since $E_k$ converges to $E$,  when $k$ is large, $E_k$ is very closed to $E$, and hence $F$ (since $E\bs B=F\bs B$) on the sphere $\partial B$. Hence near the sphere, we can use Federer Fleming projections to weld $E_k$ and $F$ together, without adding much measure. The new obtained sets are called $F_k$. They will coincide with $E_k$ outside a slightly larger ball $B'$. See Section 3 for the construction.

But the key is to prove that $F_k$ are competitors for $E_k$. For this purpose we have to proved that the homology group of $U\bs F_k$ is controlled by that of $U\bs E_k$, using the fact that the homology group of their limit $U\bs F$ is controlled by the limit of $U\bs E_k$. By some standard argument, we can restrict ourselves to only look at the homology group on the sphere $\partial B'\bs F_k$. In order to pass to the limit, we need some finiteness of homology groups. But there is no reason why the homology group has some finiteness property. So we are obliged to add the assumption that the coefficient group $G$ of the homology group is finitely generated. 

Also, no matter how close are $F_k$ and $E$ on the sphere, there may exists local elements in $H_{n-d-1}(U\bs E_k)$ that do not exist in $H_{n-d-1}(U\bs E)$, so that we cannot control them using  $E$. But we can kill these local element by adding $d-$dimensional grids to a neighborhood of $E$ and $E_k$. See Section 4 for detail.

\bigskip

\textbf{Acknowledgement:} The research leading to these results has received fundings from the European Research Council under the European UnionÕs Seventh Framework Programme (FP/2007-2013) / ERC Grant Agreement n. 291497. 

\bigskip

\section{Preliminaries}

In this section we will give necessary definitions and preliminaries.

\subsection{Basic notation and definition}

$B(x,r)$ is the open Euclidean ball with radius $r$ and centered at $x$;

$\overline B(x,r)$ is the closed ball with radius $r$ and center $x$;

$\H^d$ is the Hausdorff measure of dimension $d$;

For any two points $a,b\in \R^n$, $R_{ab}$ denotes the half line issued from the point $a$ and passing through the point $b$;

For a set $A\subset\R^n$, $A^\circ$ denotes its interior.
%

\begin{defn}[Local variant of the Hausdorff distance] For any compact set $K\subset \R^n$, and any two subsets $E$ and $F$ of $\R^n$, the local variant of the Hausdorff distance in $K$ between $E$ and $F$ is defined by
\be d_K(E,F)=\sup\{\mbox{dist}(x,F),x\in E\cap K\}+\sup\{\mbox{dist}(x,E),x\in F\cap K\}.\ee
Note that $d_K$ is not a distance.\end{defn}

\begin{defn}[Limit of closed sets]
Let $U\subset\R^n$, and $\{E_k\}_{k\in \N}$ and $E$ are closed subsets of $\R^n$. We say that $\{E_k\}_{k\in \N}$ converges to $E$ in $U$, if for each compact set $K\subset J$, 
\be \lim_{k\to\infty}d_K(E,E_k)=0.\ee
\end{defn}

\subsection{Definition of minimal sets}

In this part we introduce the general notion of minimal sets. 

\begin{defn}[Minimal sets]Let $0<d<n$ be integers, $U$ be an open set of $\R^n$, and $\cal F$ be a class of relatively closed sets in $U$. A set $E$ in $U$ is said to be minimal of dimension $d$ in $U$ with respect to the class $\cal F$ if 
\be \H^d(E\cap B)<\infty\mbox{ for every compact ball }B\subset U,\ee
and
\be \H^d(E\bs F)\le \H^d(F\bs E)\ee
for all set $F\in\cal F$.
\end{defn}

In the above definition, we usually call $\F$ a class of competitors, and the sets in $\cal F$ are called competitors. Different choices of competitor classes give (in general) different definitions of minimal sets. Note that if two competitor classes $\F_1$ and $\F_2$ satisfies $\F_1\subset\F_2$, then the minimality with respect to $\F_2$ implies immediately the minimality (of the same dimension) with respect to $\F_1$. Hence smaller competitor class gives a weaker notion of minimality. The Almgren competitor class below is somehow weak.

\begin{defn}[Almgren competitor] Let $U\subset\R^n$ be open. Let $E\subset\R^n$ be a closed set, and $d\le n-1$ be an integer. An Almgren competitor for $E$ in $U$ is a closed set $F\subset \R^n$ such that $E\bs U=F\bs U$, and that $F\cap U=\varphi_1(E\cap U)$, where $\varphi_t:U\to U$ is a family of continuous mappings such that 
\be \varphi_0(x)=x\mbox{ for }x\in U;\ee
\be\mbox{ the mapping }(t,x)\to\varphi_t(x)\mbox{ of }[0,1]\times U\mbox{ to }U\mbox{ is continuous;}\ee
\be\varphi_1\mbox{ is Lipschitz,}\ee
  and if we set $W_t=\{x\in U\ ;\ \varphi_t(x)\ne x\}$ and $\widehat W=\bigcup_{t\in[0.1]}[W_t\cup\varphi_t(W_t)]$,
then
\be \widehat W\mbox{ is relatively compact in }U.\ee
 
Such a $\varphi_1$ is called a deformation in $U$, and $F$ is also called a deformation of $E$ in $U$.
\end{defn}

Note that if $V\subset U$ are two open sets, then an Almgren competitor (or a deformation) of $E$ in $V$ is automatically an Almgren competitor (or deformation) of $E$ in $U$.

\begin{defn}[Almgren minimal sets]
Let $0<d<n$ be integers, $U$ an open set of $\R^n$. A relatively closed set $E\subset U$ is said to be Almgren minimal of dimension $d$ in $U$ if 
it is $d$-dimensional minimal with respect to the class of all Almgren competitors $F$ for $E$ in $U$.
\end{defn}
%

It can be seen that the notion of Almgren minimality involves the concept of deformation, which comes naturally from the physical intuition on the formation of soap films. Hence many people prefer this notion due to the physical background. Besides, since it is relatively weak, any regularity property for Almgren minimal sets holds also for other stronger types of minimal sets. However, since deformation is not always easy to control, we often have to prove the Almgren minimality by proving another stronger type of minimality, which is up to now the case for most minimal cones we know. For instance, the method of paired calibrations, introduced by \cite{Br91} and \cite{LM94}, is quite often used to prove minimality among a class of competitors satisfying some separation condition, called Mumford-Shah competitor. These are competitors only for codimension 1 sets. As its generalization to higher codimensions, the definition of topological competitors is the following.

\begin{defn}Let $0<d<n$ be integers. Let $G$ be an abelian group. Let $U\subset \R^n$ be an open set, $E$ be a closed set in $U$. A closed set $F\subset U$ is said to be a $d-$dimensional $G$-topological competitor for $E$ in $U$ if there exists an open set $V\subset U$, such that
\be E\bs V=F\bs V,\ee
and  for each $n-d-1$-simplical $G$-cycle in $U\bs (V\cup E)$, if it represents a non zero element in the homology group $H_{n-d-1}(U\bs E; G)$, then it also represents a non zero element in $H_{n-d-1}(U\bs F; G)$.

When the domain $U$ is fixed, we also call $F$ a $G$-topological competitor for $E$ in $V$. 
\end{defn}

\begin{rem}we are not going to say precisely which type of homology we are using,  because in our setting, the topological spaces are always very nice (open subset of $\R^n$, or the support of a simplicial complex). However in the proofs, we often use the simplicial chain for convenience.
\end{rem}

\begin{rem}As before, one can easily check that if $V_1\subset V_2\subset U$, then $F$ is a topological competitor for $E$ in $V_1$ implies that $F$ is a topological competitor for $E$ in $V_2$.\end{rem}

\begin{defn}[$G$-Topological minimal sets]\label{min}
Let $0<d<n$ be integers. Let $G$ be an abelian group. Let $U\subset \R^n$ be an open set. A closed set $E$ in $U$ is said to be $G$-topologically minimal of dimension $d$ in $U$ if it is minimal of dimension $d$ with respect to the class of all $G$-topological competitors of dimension $d$ for $E$ in balls.\end{defn}

A first relation between the two kinds of minimal sets is due to the following:

\begin{pro}Let $V\subset U\subset \R^n$ be open sets, and $\overline V\subset U$. Let $E\subset U$ be closed. Then for any coefficient group $G$, and any open set $V'\supset \overline V$, any deformation of $E$ in $V$ is a $G$-topological competitor for $E$ in $V'$. Thus the class of $G$-topological competitors of dimension $d$ for $E$ in $U$ is larger than the class of Almgren competitors of $E$ in $U$. And hence any $G$-topological minimal set of dimension $d$ in $U$ is Almgren minimal of dimension $d$ in $U$.
\end{pro}

The proof is standard, using mainly transversality. See for example the proof in Proposition 3.7 of \cite{topo}.

\subsection{Regularity of minimal sets}

In this part we cite some regularity results for reduced minimal sets that will be useful later. Some of these results were proved by many people in many ways, but for convenience the author will cite G.David's work systematically. Also, these results are proved for Almgren minimal sets. But due to Proposition 2.19, they also hold for topological minimal sets.

\begin{defn}[Reduced set] Let $U\subset \R^n$ be an open set. For every closed subset $E$ of $U$, denote by
\be E^*=\{x\in E\ ;\ H^d(E\cap B(x,r))>0\mbox{ for all }r>0\}\ee
 the closed support (in $U$) of the restriction of $H^d$ to $E$. We say that $E$ is reduced if $E=E^*$.
\end{defn}

\begin{rem}It is easy to see that
\be H^d(E\bs E^*)=0.\ee
And it is not hard to prove that a set $E$ is Almgren or topologically minimal if and only if $E^*$ is. As a result it is enough to study reduced minimal sets.
\end{rem}

\begin{thm}[Uniform Ahlfors regularity for minimal sets. See \cite{DS00} Proposition 4.1]For any pair of integers $d<n$, there exists a constant $C=C(n,d)>1$, such that the following holds: 
Let $U\subset \R^n$ be open, let $E$ be a reduced Almgren minimal set in $U$. Then for any ball $B(x,r)$ such that $x\in E$ and $B(x,2r)\subset U$, we have
\be C^{-1}r^d<\H^d(E\cap B(x,r))< Cr^d.\ee
\end{thm}

%

For proving existence for minimizers in various settings, we always need the lower semi continuity of Hausdorff measure with respect to the Hausdorff distance, that is, for a sequence of sets $E_k$ in a domain $U$ that converges (locally) to a set $E$ with respect to the Hausdorff distance, we want to have
\be \H^d(E)\le\liminf_{k\to\infty}\H^d(E_k).\ee

This does not hold in general. But if $E_k$ are reduced minimal sets, then this is true. 

\begin{thm}[cf. \cite{GD03} Theorem 3.4] Let $\Omega\subset\R^n$, $0<d<n$. Suppose that for each $k\ge 0$, $E_k$ is a reduced minimal set of dimension $d$ in $\Omega$, and that $E_k$ converges to $E$. Then
\be \H^d(E\cap W)\le\liminf_{k\to\infty}\H^d(E_k\cap W)\ee
for every open set $W\subset \Omega$.
\end{thm}

%
%
%
%
%
%
%
%
%

\subsection{Federer-Fleming Projection on dyadic complexes}

In this part we give the notations and conventions of dyadic complexes, and recall the definition of Federer-Fleming Projection on dyadic complexes. These will be used in the construction of topological competitors. The whole procedure is a typical technique in geometric measure theory.

Fix any $n\in \N$. Let $m\in \N$ be a positive integer. Denote by $\Delta_m$ the set of all (closed) dyadic cubes in $\R^n$, For $k<n$, denote by $\Delta_{k,m}$ the set of all $k$-dimensional faces of cubes in $\Delta_m$. An element in $\Delta_{k,m}$ is called a $k-$dimensional dyadic cube of length $2^{-m}$.

\begin{defn}[Dyadic complex] Let $0<k\le n$ be integers. 
Let $Q=\{\sigma_1,\sigma_2,\cdots,\sigma_l\}$ be a finite family of $k$-dimensional dyadic cubes of $\R^n$. For $0\le i\le k$ and $\sigma\in Q$ denote by $\K_i(\sigma)$ the set of all $i$-dimensional faces of the cube $\sigma$. Set $\K_i(Q)=\cup_{\sigma\in Q}\K_i(\sigma)$. Note that $\K_k(Q)=Q$.

Set $\K(Q)=\cup_{i=0}^k\K_i(Q)$.

We say that a family $\K$ of dyadic cubes of dimension at most $k$ is a dyadic complex of dimension $k$ if there exists a a subfamily $Q$ of dyadic cubes in $\Delta_{k,m}$ for some fixed $m\in \N$ such that $\K=\K(Q)$.

Obviously, if $\K$ is a dyadic complex, then

\be\forall \a,\beta\in \K, \a\ne\beta\Rightarrow \a^\circ\cap\beta^\circ=\emptyset.\ee

For $0\le i\le k$, denote by $\K_i$ the set of all its $i-$dimensional faces, and $\K^i=\cup_{j=0}^i \K_j$ the $i-$dimensional sub-complex of $\K$.
Denote by $|\K|$ the support of the complex $\K$:
\be |\K|=\bigcup_{\sigma\in \K}\sigma.\ee
The support $|\K^i|$ of the $i-$dimensional sub-complex $\K^i$ is called the $i$-skeleton of $\K$.
\end{defn}

%
%
%
Now we want to see how to project a given closed set onto faces of cubes.

\begin{defn}[Radial projection] Let $\sigma$ be a $k-$dimensional cube in $\R^n$, and $x\in \sigma^\circ$. Define the radial projection $\Pi_{\sigma,x}$ on the faces of $\sigma$ as follows:
\be\Pi_{\sigma, x}:=\left\{\begin{array}{rcl}\sigma\bs\{x\}&\rightarrow&\partial\sigma;\\
y&\mapsto&z\in R_{x,y}\cap\partial\sigma,\end{array}\right.\ee
where $R_{xy}$ denotes the half line issued from $x$ and passing through $y$.
\end{defn}

\begin{rem}Any radial projection on the faces of $\sigma$ fixes the points of $\partial \sigma$.
\end{rem}

Any radial projection $\Pi_{\sigma,x}$ is continuous on $\sigma\bs \{x\}$, and is Lipschitz on $\sigma\bs B(x,r)$ for any $r>0$. However the Lipschitz constant will blow up when $r\to 0$. Hence given a closed set $E$ contained in $\sigma$, a radial projection can enlarge the measure quite a lot. However, the following Lemma says that if we are allowed to choose the projection center, then the measure of the projection will be less than a uniform multiple of the measure of the original set.

\begin{lem}[cf.\cite{DS00} Lemma 3.22] Let $1\le d<k\le n$ be integers. There exists a constant $K=K(d,k)>0$ that only depends on $d$ and $k$, such that for any $k-$dimensional cube $\sigma\in \R^n$, and any set $E\subset \sigma$ with locally finite $d-$dimensional Hausdorff measure in $\sigma$, we can find a subset $X$ of $\sigma^\circ\bs E$ with non zero $\H^k$ measure, such that
\be \forall x\in X, \H^d(\Pi_{\sigma,x}(E))\le K\H^d(E).\ee
\end{lem}

\begin{rem}If $E$ is closed (and hence compact, because $\sigma$ is compact), then for any $x\in X\subset E^C$, the projection $\Pi_{\sigma,x}$ is Lipschitz on $E$ (but the Lipschitz constant could be very large).
\end{rem}

Let us continue on Federer-Fleming projection. By Lemma 2.35, for $d<k\le n$, for each $k-$dimensional dyadic complex $\K$, if $E\subset|\K|$ is a closed set with locally finite $d-$dimensional Hausdorff measure, then for each $k-$dimensional face $\sigma\in \K_k$, there exists a radial projection $\Pi_\sigma$ on faces of $\sigma$ such that
\be \H^d(\Pi_\sigma(E\cap\sigma))\le K(d,k)\H^d(E\cap\sigma).\ee
Then we can define $\phi_{k-1}:E\to |\K^{k-1}|$, such that
\be \phi_{k-1}|_{\sigma}=\Pi_\sigma\mbox{ for all }\sigma\in \K_k.\ee
$\phi_{k-1}$ is well defined, because when two cubes $\a,\beta$ of the same dimension meet each other, (2.30) says that they can only meet each other at their boundaries. But  $\Pi_\a$ and $\Pi_\beta$ are both equal to the identity on boundaries, hence they agree on $\a\cap\beta$.

Set $E_{k-1}=\phi_{k-1}(E)\subset|\K^{k-1}|$. Then by (2.38) we have
\be \H^d(\phi_{k-1}(E))\le K(d,k)\H^d(E).\ee

Now if $d=k-1$ we stop; otherwise in the $k-1$-dimensional complex $\K^{k-1}$, we can do the same thing for the $d-$dimensional subset $\phi_{k-1}(E)$ of $|\K^{k-1}|$, with a Lipschitz map $\phi_{k-2}:\phi_{k-1}(E)\to |\K^{k-2}|$ such that
\be \H^d(\phi_{k-2}\circ\phi_{k-1}(E))\le K(d,k)K(d,k-1)\H^d(E).\ee

We carry on this process until the map $\phi_d:\phi_{d+1}\circ\cdots\circ\phi_{k-1}(E)\to |\K^d|$ is defined, with
\be \H^d(\phi_d\circ\cdots\phi_{k-2}\circ\phi_{k-1}(E))\le K(d,k)K(d,k-1)\cdots K(d,d+1)\H^d(E).\ee

Set $\phi'=\phi_d\circ\cdots\phi_{k-2}\circ\phi_{k-1}: E\to |\K^d|$. It is Lipschitz, and $\phi'|_{\K^d}=Id$. Set $K_1(d,k)=K(d,k)K(d,k-1)\cdots K(d,d+1)$. Then we have
\be \H^d(\phi'(E))\le K_1(d,k)\H^d(E).\ee

Such a $\phi'$ is called a \textbf{radial projection} (for $d-$dimensional sets) on a dyadic complex.
  
  But we do not stop here. We want to construct a Lipschitz map $\phi:E\to|\K^d|$, such that modulo $\H^d$-null sets, the image $\phi(E)$ is a union of $d-$faces of $\K$. That is, if $\sigma\in \K_d$, then 
  \be \sigma^\circ\cap \phi(E)\ne\emptyset\Rightarrow \sigma\subset\phi(E).\ee
 Here for our map $\phi'$, the image $\phi'(E)$ may meet the interior of a $d$ face $\sigma$ of $\K$ but not contain it. To deal with this issue, for each $\sigma\in \K_d$ that does not satisfy (2.44) with the set $\phi'(E)$, take $x\in \sigma^\circ\bs \phi'(E)$, and denote by $\Pi_\sigma=\Pi_{\sigma,x}$. Then $\Pi_\sigma$ is Lipschitz on $\phi'(E)\cap\sigma$ (since $\phi'(E)$ is compact), and it sends $\phi'(E)$ to the boundary of $\sigma$, which is of dimension $d-1$. In other words, when $\phi'(E)$ does not cover the whole $\sigma$, we "clean" it out of $\sigma$ with $\Pi_\sigma$. 
 
 Define $\phi'':\phi'(E)\to|\K^d|$ as the following: for $\sigma\in\K_d$ that satisfies (2.44) with the set $\phi'(E)$, $\phi''|_\sigma=Id$, and for $\sigma\in \K_d$ that does not satisfy (2.44) with the set $\phi'(E)$, set $\phi''|_\sigma=\Pi_\sigma$. Then $\phi'':\phi'(E)\to |\K^d|$ satisfies
 \be \H^d(\phi''(\phi'(E))\le \H^d(\phi'(E)). \ee
 
 Such a $\phi''$ is called a \textbf{polyhedral erosion}.
 
Now set $\phi=\phi''\circ\phi'$. Then $\phi$ is a Lipschitz map from $E$ to $|\K^d|$ that satisfies (2.44), and
 \be \H^d(\phi(E))\le K_1(d,k)\H^d(E).\ee
 
 Such a projection $\phi$ is a \textbf{Federer-Fleming projection} for a set $E\subset|\K|$. Of course, by extension of Lipschitz functions, we have the following
 
\begin{lem}[Federer-Fleming projection]Let $1\le d<k\le n$ be integers, then there exists a constant $K_1(d,k)$ that only depends on $d$ and $k$, such that the following is true: If $\K$ is a $k-$dimensional dyadic complex, and $E\subset|\K|$ is a closed set with locally finite $d-$dimensional Hausdorff measure, then there exists Lipschitz maps $\phi',\phi''$ and $\phi$ from $|\K|$ to $|\K|$ such that

 $1^\circ$ $\phi':E\to |\K^d|$ is a radial projection, $\phi'|_{\K^d}=Id$, and satisfies (2.43);
 
 $2^\circ$ $\phi'':\phi'(E)\to |\K^d|$ is a polyhedral erosion, hence does not increase Hausdorff measure;
 
 $3^\circ$ $\phi=\phi''\circ\phi':E\to |\K^d|$ is a Federer Fleming projection that satisfies (2.44) and (2.46).
\end{lem}

In our construction, we will only deform our sets locally. That is, we will have an $n$ dimensional dyadic complex $\K$, and a set $E$ of dimension $d$ that is not contained in $|\K|$, and we want to deform that part $E\cap |\K|$ inside $\K$, while keeping $E\bs |\K|$ fixed. Notice that in this case, points on $\partial |\K|$ should be fixed as well. For that purpose, we first use $\phi_{n-1}$ on $\K$ to deform $E\cap |\K|$ to $|\K_{n-1}|$. Next, we only do the Federer-Fleming Projection on $n-1$-faces that are not on the boundary of $|\K|$. 

More precisely, for $d\le k\le n-1$, let $\K^*_k$ be the set of all $k$-faces $\sigma$ of $\K$ such that $\sigma^\circ\cap \partial |\K|=\emptyset$.

Let $\phi^*_{n-1}=\phi_{n-1}$ on $|\K|$, and $\phi^*_{n-1}=id$ outside $|\K|$. 

Now if $\phi^*_k$ is already defined, then define $\phi^*_{k-1}$ as follows:

$1^\circ$ $\phi^*_{k-1}$ is the radial projection from $\phi^*_k\circ \cdots\circ \phi^*_{n-1}(E)\cap \K^*_k$ to $\K_{k-1}$ (not necessarily in $\K^*_{k-1}$);

$2^\circ$ $\phi^*_{k-1}=Id$ on $\R^n\bs |\K|^\circ$.

We can define $\phi^*_k$ until $k=d$. Let $\phi'^*=\phi^*_d\circ \cdots\circ \phi^*_{n-1}$. As a last step, let $\phi''^*$ be the polyhedral erosion from $\phi'^*(E)\cap \K^*_d$ to $\K_{d-1}$. Let $\phi^*=\phi''^*\circ \phi'^*$. Then we have the following

\begin{lem}[local Federer-Fleming projection]Let $1\le d\le n$ be integers, then there exists a constant $K_1(d,n)$ that only depends on $d$, such that the following is true: If $\K$ is a $n-$dimensional dyadic complex in $\R^n$, and $E$ is a closed set with locally finite $d-$dimensional Hausdorff measure, then there exists Lipschitz maps $\phi'^*,\phi''^*$ and $\phi^*$ from $\R^n$ to $\R^n$ such that

$1^\circ$ $\phi'^*|_{(\R^n\bs |\K|^\circ)\cup |\K_d|}=id$, $\phi'^*(E\cap |\K|)\subset |\K^d|\cup \partial |\K|$, and 
\be \H^d(\phi'^*(E\cap |\K|))\le K_1(d,n)\H^d(E\cap |\K|).\ee
 
 $2^\circ$ $\phi''^*|_{(\R^n\bs |\K|^\circ)}=id$. The restriction $\phi''^*:{\phi'^*(E)\cup |\K^*_d|}\to |\K^*_d|$ is a polyhedral erosion, hence does not increase Hausdorff measure of $\phi'^*(E)$;
 
 $3^\circ$ $\phi^*=\phi''^*\circ\phi'^*:\R^n\to\R^n$ is a Federer Fleming projection from $E$ inside $|\K|$ that satisfies 
 \be \phi^*|_{(\R^n\bs |\K|^\circ)\cup |\K_{d-1}|}=id;\ee
 \be \phi^*(E\cap |\K|)\subset |\K_d|\cup \partial |\K|;\ee
  \be \sigma^\circ\cap \phi(E)\ne\emptyset\Rightarrow \sigma\subset\phi(E) \mbox{ for all }\sigma\in \K^*_d;\ee
  and
   \be \H^d(\phi(E))\le K_1(d,n)\H^d(E).\ee
\end{lem}

\section{The construction of competitors}

After all the preparation, we will begin to prove the main theorem.

\begin{thm}Let $G$ be a finitely generated abelian group. Let $1\le d<n$ be integers. Let $E_k$ be a sequence of reduced $d-$dimensional $G$-topological minimal sets in $U\subset\R^n$, and $E_k$ converge (in the sense of Definition 2.3) to a set $E$. Then $E$ is a reduced $d-$dimensional $G$-topological minimal set in $U$.
\end{thm}

\nd

We fix the group $G$, and topological minimal set means $G-$topological minimal sets in the whole proof.

Let $E_k,k\in \N$ be a sequence of reduced topological minimal sets of dimension $d$ in $U\subset\R^n$, and $E_k$ converge to a set $E$ . It is clear that $E$ is a reduced set as well. We want to prove that $E$ is also topologically minimal of dimension $d$. 

Suppose not. That is, there exists a ball $B_1$ with $\overline B_1\subset U$, and $F$ a topological competitor for $E$ in $B_1$, such that
\be A=\H^d(E\cap B_1)-\H^d(F\cap B_1)>0.\ee

We are going to use this set $F$ to construct sets $F_k, k\in \N$, such that for $k$ large, $F_k$ will be a better topological competitor for $E_k$, which will contradict our hypothesis that $F_k$ being topologically minimal. 

Without loss of generality, we can suppose that $B_1=B(0,r_1)$. Then since $\overline B_1\subset U$, there exists $r_2\in (r_1, r_1+\frac{1}{10}r_1)$ such that $\overline B(0,r_2)\subset U$. Set $B_2=B(0,r_2)$.

Let $m_0\in\N$ be such that $2^{-m_0}<(r_2-r_1)/100$. Denote by $Q_{m_0}$ the set of all closed dyadic cubes of length $2^{-m_0}$ that are contained in $B(0, r_1+\frac12(r_2-r_1))$, and denote by $D$ the union of all cubes in $Q_{m_0}$. Then
\be \overline B_1\subset D^\circ \subset D\subset B_2,\ee
and there exists $\e_1>0$ such that for any $r\in[1-2\e_1, 1+2\e_1]$, we also have
\be \overline B_1\subset rD^\circ\subset rD\subset B_2,\ee
where $rD=\{rx; x\in D\}$ for $r\in\R$. Moreover, $\partial D$ is a finite union of dyadic $n-1$ cubes.

\begin{lem}For any $x\in D$, and any $t<1$, $tx\in D^\circ$. \end{lem}

\nd Let $x\in D$. Then there exists a dyadic cube $\sigma\in Q_{m_0}$, such that $x\in \sigma$. Since $\sigma$ is a dyadic cube of length $2^{-m_0}$, there exists $l_1,\cdots, l_n\in \Z$, such that $\sigma=\Pi_{i=1}^n[2^{-m_0}l_i,2^{-m_0}(l_i+1)]$. Note that for any $i$, $|2^{-m_0}l_i|\ne |2^{-m_0}(l_i+1)|$. So set $a_i=\left\{\begin{array}{cc}2^{-m_0}l_i,\ &if\ |2^{-m_0}l_i|> |2^{-m_0}(l_i+1)|\\ 
2^{-m_0}(l_i+1),\ &if\ |2^{-m_0}l_i|< |2^{-m_0}(l_i+1)|\end{array}.\right.$ Then $(a_1,\cdots, a_n)$ is the (unique) farthest point in $\sigma$ from the origin. Denote by $R_\sigma$ the hyper rectangle $\Pi_{i=1}^n[-|a_i|,|a_i|]$. Then it is a union of dyadic cubes of length $2^{-m_0}$, and each of these dyadic cubes is contained in $B(0, r_1+\frac12(r_2-r_1))$, since for any $y\in R_\sigma$, $|y|\le |(a_1,\cdots, a_n)|$. By definition of the region $D$, each of these cubes are contained in $D$. Hence $R_\sigma\subset D$.

Now let $t<1$, then it is clear that for any $1\le i\le n$, its $i$-th coordinate $(tx)_i$ are such that $|(tx)_i|=t|x_i|<|x_i|\le |a_i|$, where $x_i$ denotes the $i$-th coordinate of $x$. As a result, $tx\in R_\sigma^\circ$, and hence $tx\in D^\circ$.\qed


Now let $f$ be the map $(1+\e_1)D\bs (1-\e_1)D^\circ \to [1-\e_1,1+\e_1]$, $f(x)=\inf\{r: x\in rD\}$. By Lemma 3.5, $f$ is well defined, and $f(x)=r$ if and only if $x\in\partial (rD)$.

\begin{lem} The map $f$ is $2^{m_0}$-Lipschitz.\end{lem}

\nd 

Let $x,y\in (1+\e_1)D\bs (1-\e_1)D^\circ$. Denote by $x_i$ and $y_i$ the $i$-th coordinates of $x$ and $y$ respectively.

Suppose that $f(x)=r$, that is, $x\in \partial (rD)$. Let $\sigma\in Q_{m_0}$ be such that $x\in r\sigma$. Then $x\in \partial(r\sigma)$. Define $l_1,\cdots, l_n\in \Z$, $a_i$, $R_\sigma$ as in the proof of Lemma 3.5. By definition of $a_i, |a_i|\ge 2^{-m_0}$ for each $1\le i\le n$.


Since $x\in \partial (r\sigma)$, for each $1\le i\le n$, $|x_i|/r\le |a_i|$. But $|y_i-x_i|\le  d(x,y),\forall 1\le i\le n$, hence $|y_i|\le |x_i|+d(x,y)\le r|a_i|+d(x,y)$. Therefore for each $1\le i\le n$,
\be \frac{|y_i|}{r+\frac{d(x,y)}{|a_i|}}\le|a_i|.\ee
That is, $\frac{y}{r+\frac{d(x,y)}{|a_i|}}\in R_\sigma$. Hence 
\be f(y)\le r+\frac{d(x,y)}{|a_i|}=f(x)+\frac{d(x,y)}{|a_i|}\le f(x)+2^{m_0}d(x,y)\ee
since $ |a_i|\ge 2^{-m_0}$ for each $i$. That is,
\be f(y)-f(x)\le 2^{m_0}d(x,y).\ee
By symmetry, we also have $f(x)-f(y)\le 2^{m_0}d(x,y)$. Hence $|f(x)-f(y)|\le 2^{m_0}d(x,y)$. So $f$ is $2^{m_0}$-Lipschitz.\qed

Now we apply \cite{Fe} 2.10.25 to the $2^{m_0}$-Lipschitz map $f$, and get
\be \int_{1-\e_1}^{1+\e_1}\H^{d-1}(E\cap f^{-1}(y))dy\le C2^{m_0}\H^d(E\cap(1+\e_1)D\bs (1-\e_1)D^\circ)<\infty,\ee
Hence there exists $r_0\in (1-\e_1,1+\e_1)$ such that $\H^{d-1}(E\cap f^{-1}(r))<\infty$, that is
\be \H^{d-1}(E\cap\partial (r_0D))<\infty.\ee

Without loss of generality, we can suppose that $r_0=1$. (Otherwise we can replace the sets $E$, $E_k$, $B_1$ and $B_2$ by $E/r$, $E_k/r$, $B_1/r$ and $B_2/r$, and notice that the minimality is invariant under dilatations).

Set $D_0=D$. Take $\e_2\in (0,\e_1)$ such that there exists $m_2\in \N$ such that $2^{m_2}\e_2\in \N$. In other words, $\e_2$ is a dyadic fraction. Set $D_1=(r_0-\e_2)D$, $D_2=(r_0+\e_2)D$. Then we have
\be \overline B_1\subset D_1^\circ\subset D_1\subset D_0^\circ\subset D_0\subset D_2^\circ\subset D_2\subset B_2.\ee

Moreover, for each $m\ge m_0+m_2$, the $D_i,i=0,1,2$ are all unions of dyadic cubes of length $2^{-m}$. 

Denote by $Q_m$ the set of all dyadic cubes of length $2^{-m}$ that are contained in $D_2$, $m\ge m_0+m_2$. Set $\K_m=\K(Q_m)$ (See Definition 2.29).


For $t\in (0,10^{-2}\e_1)$ (to be chosen later), denote by $Q_{m,t}$ the set of all $n-$dimensional cubes $\sigma\in Q_m$ such that there exists an $n-$dimensional cube $\sigma'\in Q_m$ such that $\sigma\cap \sigma'\ne\emptyset$, and $\sigma'\cap (E\cap (1+t)D_0\bs (1-t)D_0^\circ)\ne \emptyset$. Denote by $S_{m,t}=\K(Q_{m,t})$ the sub complex of $\K_m$ that is composed of all faces of dyadic cubes in $Q_{m,t}$. Denote by $S^d_{m,t}=\K^d(Q_{m,t})$ the $d-$dimensional sub complex of $S_m$. 

Let $Q'_{m,t}$ be the set of all $n$-dimensional cubes $\sigma\in Q_m$ such that $\sigma\cap |S_{m,t}|\ne\emptyset$. In other words, we get $Q'_{m,t}$ by adding all adjacent cubes to $Q_{m,t}$.  Define $S'_{m,t}=\K(Q'_{m,t})$ and  $S'^d_{m,t}=\K^d(Q'_{m,t})$.

By (3.11), and the continuity of the $d-$dimensional Hausdorff measure restricted to $E$, we have
\be \lim_{t\to 0}\H^d(E\cap (1+t)D_0\bs (1-t)D_0^\circ)=0.\ee

Thus, by the uniform Ahlfors regularity for topological minimal sets $E_k$, we claim that
\begin{lem} \be\lim_{t\to 0}\sup_k\H^d(E_k\cap (1+t)D_0\bs (1-t)D_0^\circ)=0.\ee\end{lem}

\nd
For each $\d>0$, there exists $t>0$ such that $\H^d(E\cap (1+t)D_0\bs (1-t)D_0^\circ)<\d$. Since $E\cap (1+t)D_0\bs (1-t)D_0^\circ$ is compact, there exists $\{B(x_i,r_i)\}_{1\le i\le N}$ a finite family of balls (with arbitrarily small radii) that cover $E\cap (1+t)D_0\bs (1-t)D_0^\circ$, and $\sum_{1\le i\le N}|r_i|^d<2\d$. By the finiteness of the family, the union $U_0:=\cup_{1\le i\le N}B(x_i,r_i)$ is an open neighborhood of $E\cap (1+t)D_0\bs (1-t)D_0^\circ$. Since $E$ is the Hausdorff limit of $E_k$, when $k$ is large enough, $E_k\cap (1+\frac t2)D_0\bs (1-\frac t2)D_0^\circ
\subset U_0$. By the uniform Ahlfors regularity (Theorem 2.24) for $E_k$, for each $1\le i\le N$, 
\be H^d(E_k\cap B(x_i,r_i))<C r_i^d,\ee
where $C=C(n,d)$ is the uniform Ahlfors regularity constant in Theorem 2.24 that only depends on $n$ and $d$.
 As a result, 
 \be H^d(E_k\cap (1+\frac t2)D_0\bs (1-\frac t2)D_0^\circ)\le\sum_{1\le i\le N}H^d (E_k\cap B(x_i,r_i))\le \sum_{1\le i\le N}Cr_i^d\le 2C\d.\ee
 This proves our Lemma 3.14.\qed

By Lemma 3.14, there exists $t_1>0$ such that for all $t<t_1$,
\be \H^d(E\cap (1+t)D_0\bs (1-t)D_0^\circ)<\frac {A}{4M},\ee
and
\be \H^d(E_k\cap (1+t)D_0\bs (1-t)D_0^\circ)<\frac {A}{4M},\ee
where $M=K_1(n,d)$ is the constant in Lemmas 2.47 and 2.48.

On the other hand, still by (3.13), we have the following:
\begin{lem} 

\be \lim_{t\to 0}\lim_{m\to\infty}\H^d(|S'^d_{m,t}|)=0.\ee
\end{lem}

\nd Fix any $t<2^{-m_0}\e_2$ small. Set $E_t=E\cap (1+t)D_0\bs (1-t)D_0^\circ$.

Take any $m\ge m_0+m_2$, and $2^{-m}<2^{-m_0}t$. For any $\sigma\in Q_m$, denote by $\xi(\sigma)=\cup\{\sigma'\in Q_m:\sigma'\cap \sigma\ne\emptyset\}$ the union of all its neighbours. Then $\xi(\sigma)$ is also a cube (but not a dyadic one).

Let $\sigma\in Q_m$ be such that $\sigma\cap E_t\ne\emptyset$. Then there exists $x\in \sigma\cap E_t$. Therefore the ball $\overline B(x,2^{-m})\subset \xi(\sigma)$. Hence 
\be \H^d(E_{2t}\cap \xi(\sigma))\ge \H^d(E_{2t}\cap \overline B(x,2^{-m}).\ee

We claim that \be E_{2t}\cap \overline B(x,2^{-m})=E\cap \overline B(x,2^{-m}).\ee

In fact, since $x\in E_t\subset (1+t)D_0\bs (1-t)D_0^\circ$, we know that $f(x)\in [1-t,1+t]$, where $f$ is defined above Lemma 3.6. By Lemma 3.6, $f$ is $2^{m_0}$-Lipschitz, hence for all $y\in \overline B(x,2^{-m})$, $|f(y)-f(x)|\le 2^{m_0}d(x,y)\le 2^{m_0}\times 2^{-m}<t$. That is, $f(y)\in [1-2t,1+2t]$. In other words, $y\in (1+2t)D_0\bs (1-2t)D_0^\circ$. Hence $\overline B(x,2^{-m})\subset (1+2t)D_0\bs (1-2t)D_0^\circ$. As a result, $E_{2t}\cap \overline B(x,2^{-m})=E\cap (1+2t)D_0\bs (1-2t)D_0^\circ\cap \overline B(x,2^{-m})=E\cap \overline B(x,2^{-m}).$ Thus we get Claim (3.23).

Recall that the sets $E_k$ are topologically minimal in $U$, and hence are Almgren minimal in $U$ (cf. Proposition 2.19). By Theorem 2.24, they are Ahlfors regular with a uniform constant $C$. Since $x\in E=\lim_{k\to\infty}E_k$, when $k$ is large
\be \H^d(E_k\cap \overline B(x,2^{-m})\ge (2C)^{-1}2^{-dm}.\ee

We want to prove that 
\be \H^d(E\cap \overline B(x,2^{-m}))\ge C'2^{-dm}\ee 
for some $C'>0$ as well. So let $\d>0$ be such that $10\d<$dist$(\overline B(x,2^{-m}),\partial U)$. For any $\e>0$ small, we can cover $E\cap \overline B(x,2^{-m})$ by countably many balls $B(y_i,t_i),i\in I$ with radius less than $\d$, such that
\be \sum_{i=1}^\infty t_i^d\le \H^d(E\cap \overline B(x,2^{-m})+\e.\ee

By Vitali covering theorem, we can find a subfamily $J\subset I$, such that the balls $B(y_i,t_i),i\in J$ are disjoint, and $E\cap \overline B(x,2^{-m})\subset \cup_{i\in J}B(y_i,5t_i)$. By compactness of $E\cap \overline B(x,2^{-m})$, we can suppose that $J$ is finite. Hence $\cup_{i\in J}B(y_i,5t_i)$ is an open neighborhood of $ E\cap \overline B(x,2^{-m})$. Therefore when $k$ is large, $E_k\cap \overline B(x,2^{-m})\subset \cup_{i\in J}B(y_i,5t_i)$ as well. Thus we have
\be \begin{split}(2C)^{-1}2^{-dm}&\le \H^d(E_k\cap \overline B(x,2^{-m}))\le\sum_{i\in J}\H^d(E_k\cap B(y_i,5t_i))\\
&\le \sum_{i\in J}C(5t_i)^d=5^dC\sum_{i\in J}t_i^d\le 5^dC(\H^d(E\cap \overline B(x,2^{-m})+\e).\end{split}\ee
 Let $\e\to 0$, we get
 \be \H^d(E\cap \overline B(x,2^{-m}))\ge C'2^{-dm},\ee
 where $C'=\frac12C^{-2}5^{-d}$. 
 
 Combine with (3.22) and (3.23), we get
 \be \H^d(E_{2t}\cap \xi(\sigma))\ge C'2^{-dm}.\ee
 
 Note that (3.29) is true for all $\sigma\in Q_m$ such that $\sigma\cap E_t\ne \emptyset$.
 
 On the other hand, note that all the $\sigma,\sigma\in Q_m$ are essentially disjoint, hence the $\xi(\sigma),\sigma\in Q$ have uniformly finite overlap that depends only on $n$. That is, there exists $C_1$ that depends only on $n$ (but not on $m$), such that
 \be \sum_{\sigma\in Q_m}1_{\xi(\sigma)}\le C_1.\ee
 As a result,
 \be \sum_{\sigma\in Q_m,\sigma\cap E_t\ne \emptyset} \H^d(E_{2t}\cap \xi(\sigma))\le C_1\H^d(E_{2t}).\ee
Combine with (3.29), we get
\be \sum_{\sigma\in Q_m,\sigma\cap E_t\ne \emptyset}2^{-dm}\le C_2\H^d(E_{2t}),\ee
for some constant $C_2>0$ that only depends on $n$ and $d$ (but not on $m$ and $t$).

Now for any $\sigma\in Q_m$, set $T(\sigma)=\{\sigma'\in Q_m,\sigma'\cap \xi(\sigma)\ne\emptyset\}$. Then the $d$-skeleton $|\K^d(T(\sigma))|$ has measure 
\be\H^d(|\K^d(T(\sigma))|)=C_3 2^{-md},\ee
where $C_3$ is a constant that only depends on $n$ and $d$. On the other hand, by definition, 
\be Q'_{m,t}=\cup_{\sigma\in Q_m,\sigma\cap E_t\ne\emptyset}T(\sigma),\ee
hence the $d$-skeleton
\be |S'^d_{m,t}|=|\K^d(Q'_{m,t})|\subset \cup _{\sigma\in Q_m,\sigma\cap E_t\ne\emptyset}|\K^d(T(\sigma))|.\ee
As a result,
\be \begin{split}
\H^d(|S'^d_{m,t}|)&\le \sum_{\sigma\in Q_m,\sigma\cap E_t\ne\emptyset}\H^d(|\K^d(T(\sigma))|)=\sum_{\sigma\in Q_m,\sigma\cap E_t\ne\emptyset}C_32^{-md}\\
&=C_3(\sum_{\sigma\in Q_m,\sigma\cap E_t\ne\emptyset}2^{-md})\le C_3C_2\H^d(E_{2t})
\end{split}\ee
by (3.33) and (3.32).

As a result,
\be\limsup_{m\to \infty}\H^d(|S'^d_{m,t}|)\le C_3C_2\H^d(E_{2t})= C_3C_2\H^d(E\cap (1+2t)D_0\bs (1-2t)D_0^\circ),\ee
which tends to zero when $t\to 0$ by (3.13). Therefore,
\be \lim_{t\to 0}\lim_{m\to \infty}\H^d(|S'^d_{m,t}|)=0.\ee

\qed

By Lemma 3.20, there exists a $\tau\in(0,\min\{2^{-10}\e_2,\frac12 t_1\})$ , and $m_3>m_0+m_2$, with $2^{-m3}<<\tau$, such that
\be \H^d(|S'^d_{m_3,\tau}|)<\frac14 [\H^d(E\cap B_1)-\H^d(F\cap B_1)]=\frac {A}{4}.\ee

We fix this pair of $m_3, \tau$. Let $Q$ denote $Q_{m_3, \tau}$, $S=S_{m_3,\tau}$, $S^d=S^d_{m_3,\tau}$; let $Q'=Q'_{m_3, \tau}$, $S'=S'_{m_3,\tau}$, $S'^d=S'^d_{m_3,\tau}$.

Since $2^{-m_3}<<\tau<\frac 12 t_1$, $|S'|\subset (1+\frac34 t_1)D_0\bs (1-\frac34 t_1)D_0^\circ$. Therefore by (3.18) and (3.19),
\be \H^d(E\cap |S'|)<\H^d(E\cap (1+\frac34 t_1)D_0\bs (1-\frac34 t_1)D_0^\circ)<\frac {A}{4M},\ee
and 
\be \H^d(E_k\cap |S'|)<\H^d(E_k\cap (1+\frac34 t_1)D_0\bs (1-\frac34 t_1)D_0^\circ)<\frac {A}{4M}.\ee

Let $Q_B$ denote the set of polygons in $Q$ that touch the boundary of $|S|$, that is, the outside layer of $Q$. 

We claim that 
\be \mbox{for any }\sigma\in Q_B\mbox{ such that }\sigma\subset (1+\tau)D_0\bs (1-\tau)D_0^\circ, \sigma\cap E=\emptyset.\ee

 In fact, for any $\sigma\in Q$ which is contained in $(1+\tau)D_0\bs (1-\tau)D_0^\circ$, if $\sigma\cap E\ne\emptyset$, then $ \sigma\cap [E\cap (1+\tau)D_0\bs (1-\tau) D_0^\circ]\ne\emptyset$. As a result, by definition of $Q$, all cubes adjacent to $\sigma$ must also belong to $Q$. Thus $\sigma$ cannot touch the boundary of $|S|=\cup\{\sigma\in Q\}$.
 
 As a result, 
 \be E\cap (1+\frac 34\tau)D_0\bs (1-\frac 34\tau) D_0^\circ\subset |S|^\circ.\ee
 
 Therefore, since $E_k$ converges to $E$, there exists $k_1$ such that for any $k>k_1$, 
\be E_k\cap (1+\frac12\tau)D_0\bs (1-\frac12\tau)D_0^\circ\subset |S|^\circ.\ee 

That is, if we denote by $Q'_B$ the set of polygons in $Q'$ that touch the boundary of $|S'|$, then for any $\sigma\in Q'_B$, 
\be \sigma\cap E_k\cap [(1+\frac12\tau)D_0\bs (1-\frac12\tau)D_0^\circ]=\emptyset\ee
 for $k>k_1$.
 
%

Therefore, when $k>k_1$, we can use the local Federer-Fleming projection (Lemma 2.48) inside $|S'|$ to project each $E_k$ to a a subset of $|S'^d|\cup\partial |S'|$. More precisely, there exists a Lipschitz map $\varphi_k: U\to U$ such that $\varphi_k|_{|S'|^C\cup \partial |S'|\cup |S'^d|}=id$, and
\be \varphi_k(E_k\cap |S|^\circ)\subset |S^d|,\ee
\be \H^d(\varphi_k(E_k\cap |S'|))\le M \H^d(E_k\cap |S'|).\ee

Also note that when $k>k_1$, by (3.44), the part of the set $E_k$ inside $(1+\frac12\tau)D_0\bs (1-\frac12\tau)D_0^\circ$ is contained in $|S|^\circ$. Hence by (3.46), 
\be \varphi(E_k\cap [(1+\frac12\tau)D_0\bs (1-\frac12\tau)D_0^\circ])\subset |S^d|.\ee

We can also do the local Federer-Fleming projection for $F$ in $|S'|$: there exists a Lipschitz map $\psi: U\to U$ such that $\varphi_k|_{|S'|^C\cup \partial |S'|\cup |S'^d|}=id$, and
\be \psi(F\cap |S|^\circ)\subset |S^d|,\ee
\be \H^d(\psi(F\cap |S'|))\le M \H^d(F\cap |S'|).\ee

We know that the set $F$ equals $E$ outside $B_1$, hence $F\cap (1+\frac 34\tau)D_0\bs (1-\frac 34\tau) D_0^\circ=E\cap (1+\frac 34\tau)D_0\bs (1-\frac 34\tau) D_0^\circ$. By (3.43), 
\be F\cap (1+\frac 34\tau)D_0\bs (1-\frac 34\tau) D_0^\circ=E\cap (1+\frac 34\tau)D_0\bs (1-\frac 34\tau) D_0^\circ\subset |S|^\circ.\ee
thus
\be \psi(F\cap [(1+\frac34\tau)D_0\bs (1-\frac34\tau)D_0^\circ])\subset |S^d|.\ee


%
%

Now define $F_k=(\varphi_k(E_k)\bs D_0^\circ)\cup (\psi(F)\cap D_0)\cup |S'^d|.$ That is, we use $|S'^d|$ to weld the part of $\varphi_k(E_k)$ outside $D_0$, the part of $\psi(F)$ inside $D_0$ together. We can do this because by (3.46) and (3.49), on $\partial D_0$, both $\varphi_k(E_k)$ and $\psi(F)$ are contained in $|S'^d|\cap \partial D_0$. Note that by definition, 
\be F_k\bs D_2^\circ=E_k\bs D_2^\circ.\ee 

Now we estimate the measure of $F_k\cap D_2^\circ$, note that
\be\varphi_k(E_k)\bs D_0^\circ\subset [(E_k\bs |S'|)\bs D_0^\circ]\cup (\varphi_k(E_k\cap |S'|)),\ee
and by (3.41) and (3.47), 
\be \H^d(\varphi_k(E_k\cap |S'|))\le M \H^d(E_k\cap |S'|)\le M\times\frac{A}{4M}=\frac A4;\ee
Also, by (3.44), $E_k\cap (1+\frac 12\tau)D_0\bs (1-\frac 12\tau)D_0^\circ\subset |S'|$, 
hence $(E_k\bs |S'|)\bs D_0^\circ\subset U\bs (1+\frac 12\tau)D_0$,
therefore \be \H^d ([\varphi_k(E_k)\bs D_0^\circ]\cap D_2^\circ)\le \H^d((E_k\cap D_2^\circ\bs(1+\frac 12\tau)D_0)+\frac A4;\ee

On the other hand, 
\be \psi(F)\cap D_0=(F\cap D_0\bs |S'|)\cup \psi (F\cap |S'|).\ee
Note that $F\cap |S'|=E\cap |S'|$, hence by (3.40) and (3.50)
\be \H^d(\psi (F\cap |S'|))\le M\H^d(F\cap |S'|)= M\H^d(E\cap |S'|)\le M\times\frac{A}{4M}=\frac A4.\ee
Also, by (3.51),  $(F\bs |S'|)\cap D_0\subset F\cap (1-\frac 34\tau)D_0$,
Therefore
\be \H^d(\psi(F)\cap D_0)\le \H^d(F\cap (1-\frac 34\tau)D_0)+\frac A4.\ee

Recall that by (3.39), $\H^d(|S'^d|)\le \frac A4$. By (3.56), (3.59), and the definition of $F_k$, we have
\be \begin{split}
\H^d(F_k\cap D_2^\circ)&\le \H^d([\varphi_k(E_k)\bs D_0^\circ)+\H^d(\psi(F)\cap D_0)+\H^d(|S'^d|)\\
&\le [\H^d((E_k\cap D_2^\circ\bs (1+\frac 12\tau)D_0))+\frac A4]+[\H^d(F\cap (1-\frac 34\tau)D_0)+\frac A4]+\frac A4\\
&\le\H^d((E_k\cap D_2^\circ\bs (1+\frac 12\tau)D_0))+ \H^d(F\cap (1-\frac 34\tau)D_0)+\frac {3A}{4}.\end{split}\ee

Recall that $F$ is a competitor for $E$ in $B_1$, and $B_1\subset(1-\frac34\tau) D_0$, hence $F\cap (1-\frac34\tau) D_0 \bs B_1=E\cap (1-\frac34\tau) D_0\bs B_1 $, thus
\be  \H^d(E\cap (1-\frac34\tau) D_0)-\H^d(F\cap (1-\frac34\tau) D_0)=\H^d(E\cap B_1)-\H^d(F\cap B_1)=A.\ee

Hence (3.60) becomes
\be \H^d(F_k\cap D_2^\circ)\le \H^d((E_k\cap D_2^\circ\bs (1+\frac 12\tau)D_0))+ \H^d(E\cap (1-\frac 34\tau)D_0)-\frac A4.\ee

But $E$ is the Hausdorff limit of $E_k$, by the lower semi continuity of Hausdorff measure for minimal sets  (cf. Theorem 2.27), we have
\be \H^d(E\cap D_0)\le\liminf_{k\to\infty} \H^d(E_k\cap (1+\tau) D_0).\ee

As a result, there exists $k_2>k_1$ such that for any $k>k_2$, 
\be \H^d(E\cap (1-\frac 34\tau)D_0)\le \H^d(E_k\cap (1+\frac12 \tau) D_0)+\frac A8.\ee

Combine with (3.62), we have
\be \begin{split}H^d(F_k\cap D_2^\circ)&\le \H^d((E_k\cap D_2^\circ\bs (1+\frac 12\tau)D_0))+ \H^d(E_k\cap (1+\frac12 \tau) D_0+\frac A8-\frac {A}{4}\\
&=\H^d(E_k\cap D_2^\circ)-\frac A8\end{split}\ee
whenever $k>k_2$.

Now we have constructed the sequence $F_k$, which have smaller measure than $E_k$ when $k>k_2$. To complete the proof of Theorem 3.1, we have to prove that each $F_k$ is a topological competitor for $E_k$ in $D_2^\circ$. We will do this in the next section.

%
%

\section{$F_k$ is a topological competitor for $E_k$}

In this section, we prove the following lemma to complete the proof of Theorem 3.1.

\begin{lem} For each $k>k_2$, $F_k$ is a topological competitor for $E_k$ in $D_2^\circ$.\end{lem}

For any $k>k_2$. Denote by $E_k'=\varphi_k(E_k)\cup |S'^d|$, and $F'=\psi(F)\cup |S'^d|$. 


By definition of $F_k$, we have
\be F_k\bs D_0^\circ=E_k'\bs  D_0^\circ, F_k\cap D_0=F'\cap D_0.\ee

Recall that $\partial D_0$ is a union of $n-1-$faces of $\K_{m_3}=\K(Q_{m_3})$. So let $T$ denote the $n-1$-sub complex of $\K_{m_3}$: $T:=\{\sigma\in \K_{m_3}:\sigma\subset \partial D_0\}.$ 
 Denote by $T'$ the $n-1$-sub complex of $S'$ and $T$: $T':=\{\sigma\in S',\sigma\subset \partial D_0\}$. Then $T'=T\cap S'$, and for any $k>k_2$,
\be F_k\cap \partial D_0=E_k'\cap \partial D_0=F'\cap\partial D_0=|T'^d|.\ee

Now $\partial D_0=|T|$, and $T'^d$ is a sub complex of $T$, hence $H_0=H_{n-d-1}(\partial D_0\bs |T'^d|;G)$ is a finitely generated abelian group. 

Since $F'\cap\partial D_0=|T'^d|$, we have the natual inclusion map $j: \partial D_0\bs |T'^d|\to D_0\bs F'$, which induces a group homomorphism $j_*: H_{n-d-1}(\partial D_0\bs |T'^d|;G)\to H_{n-d-1}(D_0\bs F';G)$. Let $H=ker j_*$, then $H$ is a subgroup of $H_0=H_{n-d-1}(\partial D_0\bs |T'^d|;G)$, hence is also finitely generated. Let ${\cal A}=\{a_i,1\le i\le N\}$ be a finite set of generators of $H$.

\begin{lem}For each $a_i$, there exists a smooth simplicial $n-d-1$-cycle $\gamma_i\subset \partial D_0\bs |T'|$ (not only $|T'^d|$) which represents $a_i$. \end{lem}

\nd For the pair of topological spaces $(|T'|,|T'^d|)$, and for any $q\le d-1$, we have the exact sequence
\be H_{q+1}(|T'|,|T'^d|; \Z)\to H_q(|T'^d|; \Z)\stackrel{i_*}{\to} H_q(|T'|; \Z)\to H_q(|T'|,|T'^d|; \Z),\ee
where $i_*$ is induced by the inclusion map $i: |T'^d|\to |T'|$. 

However, for any $q\le d-1$, and any simplicial $q-$chain or $q+1$ chain is of dimension less or equal than $d$. By Federer Fleming projection, any simplicial $q-$chain or $q+1$ chain in $|T'|$ with boundary in $|T'^d|$ is homotopic to chains in $|T'^d|$. That is, any simplicial $q-$chain or $q+1$ chain in the pair $(|T'|,|T'^d|)$ represents a zero element. Hence $H_q(|T'|,|T'^d|; \Z)=H_{q+1}(|T'|,|T'^d|; \Z)=0$. As a result, the map $i_*$ in (4.5) is an isomorphism.

By the universal coefficient theorem for cohomology and the naturality, we have the following commutative diagram: 
\be\begin{array}{ccccccccc}
0       &\longrightarrow &Ext(H_{d-2}(|T'|; \Z),G) &\longrightarrow  &H^{d-1}(|T'|;G)       &\longrightarrow &Hom(H_{d-1}(|T'|; \Z),G)   &\longrightarrow &0\\
        &                   &\downarrow Ext(i_*,j)                &                   &\downarrow i^*        &                  &\downarrow Hom (i_*,j)                   &            &\\
0       &\longrightarrow &Ext(H_{d-2}(|T'^d|; \Z),G) &\longrightarrow  &H^{d-1}(|T'^d|;G)    &\longrightarrow &Hom(H_{d-1}(|T'^d|; \Z),G)&\longrightarrow &0,
\end{array}\ee
where $j$ is the identity map of $G$.

Since $i_*$ in (4.5) is an isomorphism for $d-1$ and $d-2$, the two maps $Ext(i_*,j)$ and $Hom (i_*,j)$ are isomorphisms. By the five lemma, the map 
\be H^{d-1}(|T'|;G)\stackrel{i^*}{\to}H^{d-1}(|T'^d|;G)\ee
is also an isomorphism.

Now since $\partial D_0$ is topologically an $n-1$-
sphere, by Alexander duality and its naturality with respect to inclusions, we have the commutative diagram
\be\begin{array}{ccc}
H_{n-d-1}(\partial D_0\bs |T'|;G)&\stackrel{i'_*}{\longrightarrow} &H_{n-d-1}(\partial D_0\bs |T'^d|;G)\\
\downarrow &      &\downarrow\\
H^{d-1}(|T'|;G)&\stackrel{i^*}{\longrightarrow} &H^{d-1}(|T'^d|;G),
\end{array}
\ee
where  $i': \partial D_0\bs |T'|\to \partial D_0\bs |T'^d|$ is the inclusion map. As a result, $i'$ also induces a isomorphism
\be H_{n-d-1}(\partial D_0\bs |T'|;G)\stackrel{i'_*}{\cong}H_{n-d-1}(\partial D_0\bs |T'^d|;G).\ee

That proves that each $a_i\in \A$ can be represented by a simplicial $n-d-1$-cycle $\gamma_i$ that does not touch $|T'|$. \qed

Let $V\subset D_2^\circ\bs D_1$ be a small neighborhood of $|S'|$, and $V\cap \gamma_i=\emptyset$ for any $1\le i\le N$.

Then by definition, $\varphi_k(E_k)$ (resp. $\psi(F)$) is a deformations of $E_k$ (resp. $F$) in $V$. Hence by Proposition 2.19, $\varphi_k(E_k)$ (resp. $\psi(F)$) is a topological competitor for $E_k$ (resp. $F$) in $V$. And so is $E_k'$ (resp. $F'$).

Recall that each $\gamma_i$ represents a zero element in $H_{n-d-1}(D_0\bs F'; G)$, and hence in $H_{n-d-1}(U\bs F'; G)$. But $F'$ is a topological competitor for $F$ in $V$, and $\gamma_i\in U\bs ( V\cup F)$, hence $\gamma_i$ represents also a zero element in $H_{n-d-1}(U\bs F; G)$. Recall that $F$ is a topological competitor for $E$ in $B_1$, and $\gamma_i\in D_1^C\subset B_1^C$, hence $\gamma_i$ represents also a zero element in $H_{n-d-1}(U\bs E; G)$. As a result, for each $1\le i\le N$, there exists a $n-d$-chain $\Gamma_i\subset U\bs E$ such that $\partial \Gamma_i=\gamma_i$. Denote by $|\Gamma_i|$ the support of $\Gamma_i$, then $\Gamma:=\cup_{1\le i\le N}|\Gamma_i|$ is compact, and does not touch $E$. Now since $E$ is the Hausdorff limit of $E_k$, there exists $k_3>k_2$, such that for all $k>k_3$, $E_k\cap \Gamma=\emptyset$.

Now we are ready to prove that for any $k>k_3$, $F_k$ is a $d$-dimensional topological competitor for $E_k$ in $D_2^\circ$. 

So fix any $k>k_3$. By (3.53), $F_k\bs D_2^\circ=E_k\bs D_2^\circ$. 

Let $\sigma$ be a simplicial  $n-d-1$ chain in $U\bs (D_2^\circ\cup F_k)$, and represents a zero element in $H_{n-d-1}(U\bs F_k; G)$. We want to prove that 
\be \sigma\mbox{ also represents a zero element in }H_{n-d-1}(U\bs E_k; G).\ee

Let $\Sigma$ be an $n-d$ chain in $U\bs F_k$, such that $\partial \Sigma=\sigma$. 

If $\Sigma\cap D_0=\emptyset$, then $\Sigma\subset U\bs (F_k\cup D_0)$. But by (4.2), $F_k\bs D_0=E_k'\bs D_0$, hence $\Sigma\subset U\bs (E_k'\cup D_0))\subset U\bs E_k'$, and hence $\sigma$ also represents a zero element in $H_{n-d-1}(U\bs E_k'; G)$. By Proposition 2.19, $E_k'$ is a topological competitor for $E_k$ with respect to $V$, and $\sigma\in U\bs D_2^\circ\subset U\bs V$, therefore $\sigma$ also represents a zero element in $H_{n-d-1}(U\bs E_k; G)$. 

Otherwise, $\Sigma\cap D_0\ne \emptyset$. By transversality, we can suppose that $\Sigma$ intersects $\partial D_0$ transversally. Hence the intersection $\sigma_0$ is also a simplicial $n-d-1$ cycle on $D_0\bs F_k$, and the part $\Sigma$ inside $D_0$ is a simplicial $n-d$ chain $\Sigma_0$, such that $\partial \Sigma_0=\sigma_0$. 

Since $\partial\Sigma_0=\sigma_0$, $\sigma_0$ represents a zero element in $H_{n-d-1}(D_0\bs F_k;G)$. By (4.2), $D_0\cap F_k=D_0\cap F'$, hence $\sigma_0$ represents a zero element in $H_{n-d-1}(D_0\bs F';G)$. Note that $\sigma_0\subset \partial D_0\bs F'=\partial D_0\bs |T'^d|$, hence $\sigma_0$ represents an element in the group $H$. So there exists $g_1,\cdots g_N\in G$, such that $\sigma_0$ represents the element $\sum_{1\le i\le N}g_ia_i$ in $H$. By Lemma 4.4, $\sigma$ is homologue to $\sigma_1=\sum g_i\gamma_i$ in $\partial D_0\bs |T'^d|$, and hence also in $U\bs E_k'$, since $\partial D_0\bs |T'^d|\subset U\bs E_k'$. Moreover, $\sigma_1\subset U\bs V$.

Also denote by $\Sigma_2$ the part of $\Sigma$ outside $D_0^\circ$, which is also a simplicial $n-d$ chain, and $\partial\Sigma_2=\sigma-\sigma_0.$ This means, $\sigma-\sigma_0$ represents a zero element in $H_{n-d-1}(U\bs (D_0^\circ\cup F_k); G)$. But by (4.2), $F_k\bs D_0^\circ=E_k'\bs D_0^\circ$, hence $\sigma-\sigma_0$ represents a zero element in $H_{n-d-1}(U\bs (D_0^\circ\cup E_k'); G)$, and hence in $H_{n-d-1}(U\bs E_k'; G)$. Recall that $\sigma_0$ is homologue to $\sigma_1$ in $U\bs E_k'$, hence $\sigma-\sigma_1$ represents a zero element in $H_{n-d-1}(U\bs E_k'; G)$.

Recall that $E_k'$ is a topological competitor for $E_k$ in $V$, and $\sigma-\sigma_1$ does not touch $V$,  As a result, $\sigma-\sigma_1$ represents a zero element in $H_{n-d-1}(U\bs E_k, G)$, hence there exists an $n-d$-chain $\Sigma_2'\subset U\bs E_k$ such that $\partial\Sigma_2'=\sigma-\sigma_1.$



Define the $n-d-$chain $\Theta=\sum_{1\le i\le N}g_i\Gamma_i+\Sigma_2'$. Then $\Theta$ does not touch $E_k$ since $k>k_3$, and $\partial \Theta=\sigma$. This proves (4.10). Hence the proof of Theorem 3.1 is completed. \qed

\section{Possible applications}

\subsection{Classification of singularities}

An immediate consequence of Theorem 3.1 is the following:

\begin{cor}Let $U\subset\R^n$,  and $G$ be a finitely generated group. Let $E$ be a reduced $G$-topological minimal set of dimension $d$ in $U$. Then given a point $x\in E$, any blow-up limit of $E$ at $x$ is a $G$-topological minimal cone of dimension $d$. 
\end{cor}

Here a blow-up limit of $E$ at $x$ is the limit of any converging sequence $\frac{1}{r_k}(E-x)$ with $r_k\to 0$. It describes the asymptotic behavior of $E$ around $x$ at small scales. The study of blow-up limits for sets is the key point in the classification of singularities for minimal sets.

\nd It is known that any blow-up limit of an Almgren minimal set is an Almgren minimal cone (cf. \cite{DJT} Proposition 7.31). Since topological minimal sets are all Almgren minimal, their blow-ups are cones. The corollary follows hence directly from Theorem 3.1.\qed

\subsection{Bernstein type problem}

Similarly, we can also apply Theorem 3.1 to the Bernstein type problem for minimal sets, that is, whether all topological minimal sets are cones.

The basic idea to study this problem is to look at the blow-in limits for topological minimal sets, that is, limits of the sets
\be E_r=\frac 1r E,r\to\infty.\ee

The blow-in limits for a set $E$ describe what the set $E$ looks like at infinity. And by Theorem 3.1, these blow-in limits are topological minimal cones. Then the rest of the task is to use minimality of sets to control their topological behaviors at small scales by their behaviors at large scales. See for example \cite {DJT} Section 18, \cite{globalT},\cite{global2p}, for details. 

Note that this Bernstein type problem is a typical interest for all kinds of minimizing problems in geometric measure theory and calculus of variations. One can refer to \cite{Be, Mo89, Mo86,DMS} for results on complete 2 dimensional minimal surfaces in $\R^3$, area or size minimizing currents in $\R^n$, and global minimizers for the Mumford-Shah functional.

\subsection{Local almost-Almgren minimality for product of an Almgren minimal set and $\R^n$}

Up to now, we do not know any example of Almgren minimal cone which is not topologically minimal. In fact, we only know a few Almgren minimal cones, among which the only possible non topological minimal ones are unions of almost orthogonal planes. (The author proved in \cite{2p} the Almgren minimality for a family of unions of almost orthogonal planes, and she proved then in \cite{2ptopo} that a subfamily is topologically minimal.) 

It would be of course interesting if there were any Almgren minimal cone which is not topological minimally, according to the above corollary. On the other hand, if all Almgren minimal cones are topologically minimal, things might be even better, because then many good properties for topological minimal sets could be proved in an asymptotic way for almgren minimal sets, by compactness argument and using Theorem 3.1.

Here is an example: we do not know whether the product of two Almgren minimal sets is still Almgren minimal, although this sounds reasonable. We do not even know whether the product of an Almgren minimal set with $\R$ is minimal. However the last property is true for topological minimal sets. So if all blow-up limits for Almgren minimal sets are topologically minimal, we can conclude that any blow-up limit for the product $E\times\R$ of an Almgren minimal set $E$ and $\R$ is topologically minimal. 

Of course this property alone does not guarantee anything: all manifolds admit planes (which are topologically minimal of course) as blow-up limits, but they are by no means minimal.

However for our particular example, since $E$ is Almgren minimal, $E\times \R$ admits some other useful properties. These properties could help to prove the asymptotic Almgren minimality for $E\times \R$, by a compactness argument with the help of Theorem 3.1.

\renewcommand\refname{References}
\bibliographystyle{plain}
\bibliography{reference}

\begin{thebibliography}{10}

\bibitem{Al76}
F.~J. Almgren.
\newblock Existence and regularity almost everywhere of solutions to elliptic
  variational problems with constraints.
\newblock {\em Memoirs of the American Mathematical Society}, 4(165), 1976.

\bibitem{Be}
S.~Bernstein.
\newblock Sur un th\'eor\`eme de g\'eom\'etrie et ses applications aux
  \'equations aux d\'eriv\'ees partielles du type elliptique.
\newblock {\em Comm. Soc. Math. de Khardov}, 15:38--45, 1915-17.

\bibitem{Br91}
Kenneth~A Brakke.
\newblock Minimal cones on hypercubes.
\newblock {\em Journal of Geometric analysis}, 1(4):329--338, 1991.

\bibitem{GD03}
Guy David.
\newblock Limits of {A}lmgren-quasiminimal sets.
\newblock {\em Proceedings of the conference on Harmonic Analysis, Mount
  Holyoke, A.M.S. Contemporary Mathematics series}, 320:119--145, 2003.

\bibitem{DMS}
Guy David.
\newblock {\em Singular sets of minimizers for the Mumford-Shah functional}.
\newblock Progress in Mathematics. Birkh\"auser, 2005.

\bibitem{DJT}
Guy David.
\newblock H\"older regularity of two-dimensional almost-minimal sets in $\r^n$.
\newblock {\em Annales de la Facult\'e des Sciences de Toulouse},
  XVIII(1):65--246, 2009.

\bibitem{DS00}
Guy David and Stephen Semmes.
\newblock Uniform rectifiablilty and quasiminimizing sets of arbitrary
  codimension.
\newblock {\em Memoirs of the A.M.S.}, 144(687), 2000.

\bibitem{Fe}
Herbert Federer.
\newblock {\em Geometric measure theory}.
\newblock Grundlehren der Mathematishen Wissenschaften 153. Springer Verlag,
  1969.

\bibitem{LM94}
Gary Lawlor and Frank Morgan.
\newblock Paired calibrations applied to soap films, immiscible fluids, and
  surface or networks minimizing other norms.
\newblock {\em Pacific journal of Mathematics}, 166(1):55--83, 1994.

\bibitem{global2p}
Xiangyu Liang.
\newblock Regularity for minimal sets near a union of two planes.
\newblock {\em to appear in Annales de L'Institut Fourier, arXiv: 1203.0560},
  2012.

\bibitem{2p}
Xiangyu Liang.
\newblock Almgren-minimality of unions of two almost orthogonal planes in
  $\r^4$.
\newblock {\em Proceedings of the London Mathematical Society},
  106(5):1005--1059, 2013.

\bibitem{topo}
Xiangyu Liang.
\newblock Topological minimal sets and existence results.
\newblock {\em Calculus of Variations and Partial Differential Equations},
  47(3-4):523--546, 2013.

\bibitem{YXY}
Xiangyu Liang.
\newblock Almgren and topological minimality for the set ${Y}\times {Y}$.
\newblock {\em Journal of Functional Analysis}, 266(10):6007--6054, 2014.

\bibitem{globalT}
Xiangyu Liang.
\newblock Global regularity for minimal sets near a $\mathbb{T}$ set and
  counterexamples.
\newblock {\em Revista Matem\'atica Iberoamericana}, 30(1):203--236, 2014.

\bibitem{2ptopo}
Xiangyu Liang.
\newblock On the topological minimality of unions of planes of arbitrary
  dimension.
\newblock {\em International Mathematics Research Notices 2015, doi:
  10.1093/imrn/rnv059}, 2015.

\bibitem{Mo86}
Frank Morgan.
\newblock Harnack type mass bounds and bernstein theorems for area-minimizing
  flat chains modulo $\nu$.
\newblock {\em Communications in Partial Differential Equations},
  11(12):1257--1283, 1986.

\bibitem{Mo89}
Frank Morgan.
\newblock Size-minimizing rectifiable currents.
\newblock {\em Invent.math.}, 96:333--348, 1989.

\end{thebibliography}

\end{document}